\documentclass[12pt,dvips]{amsart}
\usepackage{euler, amsfonts, amssymb, latexsym, epsfig,epic}

\newcommand\lie[1]{{\mathfrak #1}}

\newcommand\Tr{{\rm Tr\,}}
\newcommand\diag{{\rm diag}}

\newtheorem{Theorem}{Theorem} 
\newtheorem{Proposition}{Proposition} 
\newtheorem{Lemma}{Lemma}

\newtheorem*{Corollary}{Corollary}
 
\newtheorem*{Theorem*}{Theorem}
\theoremstyle{remark}

\newcommand\union{\bigcup}

\newcommand\Id{{\bf 1}}

\newcommand\complexes{{\mathbb C}}

\newcommand\naturals{{\mathbb N}}

\theoremstyle{plain}

\newtheorem*{Conjecture}{Conjecture}

\newcommand\dfn{\bf} 

\newcommand\gln{\lie{gl}_n}

\newcommand\thmlabel[1]{\label{thm:#1}}
\newcommand\proplabel[1]{\label{prop:#1}}
\newcommand\seclabel[1]{\label{sec:#1}}
\newcommand\sseclabel[1]{\label{ssec:#1}}
\newcommand\lemlabel[1]{\label{lem:#1}}

\newcommand\thmref[1]{theorem \ref{thm:#1}}
\newcommand\propref[1]{proposition \ref{prop:#1}}
\newcommand\secref[1]{section \ref{sec:#1}}
\newcommand\ssecref[1]{subsection \ref{ssec:#1}}
\newcommand\lemref[1]{lemma \ref{lem:#1}}

\newcommand\Thmref[1]{Theorem \ref{thm:#1}}

\begin{document}
\pagestyle{plain}

\title{Some schemes related to the commuting variety}
\author{Allen Knutson}
\thanks{AK was supported by the NSF and the Sloan Foundation.}
\email{allenk@math.berkeley.edu}
\date{\today}

\begin{abstract}
  The {\em commuting variety} is the pairs of $n\times n$ matrices $(X,Y)$
  such that $XY = YX$. 
  We introduce the {\dfn diagonal commutator scheme},
  $\big\{ (X,Y)  : XY-YX \hbox{ is diagonal} \big\}$,
  which we prove to be a reduced complete intersection, one component
  of which is the commuting variety. (We conjecture there to be only
  one other component.)

  The diagonal commutator scheme has a flat degeneration to the scheme
  $\big\{ (X,Y) : XY $ lower triangular, $YX$ upper triangular$\big\}$,
  which is again a reduced complete intersection, this time with 
  $n!$ components (one for each permutation).
  The degrees of these components give interesting invariants 
  of permutations.
\end{abstract}

\maketitle

\tableofcontents

\section{Statements of results}\seclabel{intro}

The {\dfn commuting variety} of a reductive Lie algebra $\lie{g}$ is
defined as the reduced subscheme of $\lie{g}\oplus \lie{g}$ cut out by
the equations $[X,Y] = 0$. It is known to be irreducible
\cite{irreducible}, but even for $\lie{g} = \gln$ it is not presently
known whether these equations serve to define it as a scheme (i.e. whether
the ideal generated by these $\dim \lie{g}$ equations is radical).

Let $\lie{h}$ be a chosen Cartan subalgebra of $\lie{g}$ (we will work
over $\complexes$ in this paper, so any two Cartan subalgebras are conjugate).
We introduce the {\dfn diagonal commutator scheme}
$$ D = \big\{ (X,Y) \in \lie{g}\oplus \lie{g} : [X,Y] \in \lie{h} \big\}. $$
This is now defined by only $\dim \lie{g} / \lie{h}$ equations,
rather than the $\dim \lie{g}$ (and who knows how many more) equations
needed to define the commuting variety.

While our first theorem could be stated for general $\lie{g}$, we only
prove it in the case $\lie{g}=\gln$, with $\lie{h}$ the diagonal matrices
(hence the name).

\begin{Theorem}\thmlabel{diagcomm}
  The diagonal commutator scheme
  $D$ is a reduced complete intersection in $\gln \oplus \gln$,
  one component of which is the commuting variety.
\end{Theorem}

It is a well-known open problem, due to Artin and Hochster, 
to show that the commuting variety
of $\gln$ is Cohen-Macaulay. (Our reference is chapter 9 of \cite{V}.)
The theorem above doesn't address that directly, but implies something related:

\begin{Corollary}
  The commuting variety is Cohen-Macaulay if and only if the union
  of the other components is Cohen-Macaulay.
\end{Corollary}

In fact we conjecture that for any reductive $\lie{g}$, the diagonal
commutator scheme has exactly two components -- i.e. only one component 
other than the commuting one.

\begin{proof}
  These two schemes are ``directly linked,'' 
  and then the theory of linkage (\cite{E}, theorem 21.23) relates
  the Cohen-Macaulayness of one component to the union of the rest.
\end{proof}

\newcommand\rhocek{\check\rho}
We will prove this by studying a certain flat
degeneration of the $\gln$ diagonal commutator scheme. 
Let $\rhocek : \complexes^\times \to GL_n(\complexes)$
denote the one-parameter subgroup
$$ \rhocek(t) = 
\left(\begin{array}{ccccc}1 & & & & \\
                          & t & & & \\
                          & & t^2 & & \\
                          & & & \ddots & \\
                          & & & & t^{n-1}\end{array}\right) $$
of diagonal matrices. 
Define $D^z$ to be the subscheme of $\gln\oplus\gln$
$$ D^z = \big\{ (X',Y') : X'Y' = \rhocek(z) Y'X' \rhocek(z^{-1}) 
\hbox{ off the diagonal}
\big\} $$
which is plainly isomorphic to $D^1 = D$ under the map
$(X',Y') \mapsto \big(\rhocek(z)\, X, \, Y \, \rhocek(z^{-1})\big)$.

The flat family $\{D^z\}, z\in\complexes^\times$ has a unique
extension to a flat family over $\complexes$, with special fiber $D^0$.
This sort of flat limit -- by rescaling some of the coordinates of 
the ambient vector space -- is called a {\em Gr\"obner degeneration}.
We know from \cite{KS} that such a limit must be again
equidimensional (up to embedded components). But in our instance much
more is true:

\begin{Theorem}\thmlabel{degen}
  The flat limit $D^0 := \lim_{z\to 0} D^z$ of this family
  is again a reduced complete intersection.

  Moreover, it is defined by the limiting equations on $D^z$,
  $$ \big\{ (X,Y) : XY \hbox{ lower triangular, } 
          YX \hbox{ upper triangular} \big\}. $$
\end{Theorem}

Ordinarily more equations are needed in such a limit; we give a
typical example. Consider the $z\to 0$ limit of the equations
$X = 0$, $X = zY$ in the $X,Y$ plane. 
For each nonzero $z$, these two equations describe two lines
intersecting at the origin.  But in the limit $z=0$ the two lines are
equal and the condition $Y=0$ is lost; the correct limit is only
obtained if that equation $Y=0$ is added to the list.
The second conclusion of the above theorem says that this unfortunate
phenomenon doesn't occur in our case: our list of
$\dim \lie{g} / \lie{h}$ equations is already enough for this limit.

(For readers familiar with Gr\"obner bases: this list of equations is
{\em not} a Gr\"obner basis, but is ``Gr\"obner enough'' for this limit
defined by a partial term order.)

We have a better handle on the components of this scheme $D^0$, because of
the large group 
$$ B_- \times B_+ = \big\{ (L,U) \in GL_n(\complexes)^2 
        : L \hbox{ lower triangular, } 
          U \hbox{ upper triangular} \big\} $$
acting on it by the rule
$$ (L,U) \cdot (X',Y') = (L \,X' \, U^{-1}, U\, Y'\, L^{-1}). $$
This group is of dimension $n^2 + n$, slightly larger than the 
$GL_n(\complexes)$ acting on the commuting variety. Since this group is
connected, it preserves (and acts on) each component of $D^0$.

For $\pi$ an $n\times n$ permutation matrix, define $D^0_\pi$ as the closure
$$ D^0_\pi 
:= \overline{ (B_-\times B_+) \cdot \union_{t,s \in H} (\pi t, s \pi^{-1}) }.$$

\begin{Theorem}\thmlabel{comps}
  The components of $D^0$ are exactly the $\{D^0_\pi\}$, $\pi\in S_n$.

  Inside the flat family $\{D^z\}_{z\in\complexes}$, consider the component
  whose generic fiber is the commuting variety as a subfamily.
  The special fiber of this subfamily is
  $D^0_{\Id}$, plus possibly some nonreduced structure.
\end{Theorem}

In the rest of the paper we prove these statements, in reverse order.
\Thmref{degen} builds on the first half of \thmref{comps}, whose proof
uses simple facts about matrix Schubert varieties (our reference is
\cite{MS}), and we prove them together in \ssecref{degencomponents}.
\Thmref{diagcomm} is then a consequence of \thmref{degen}.

We close with a number of conjectures, generalizing the standard
ones about the commuting variety.

I am pleased to thank Mark Haiman for introducing me to the
commuting variety, and many useful conversations since.
I am also grateful to David Eisenbud, Ezra Miller, and especially
Terry Tao for their insights. 

\section{A Gr\"obner degeneration using $\rhocek$}

Consider the scheme of pairs of matrices
$$ E := \big\{ (X,Y) : XY,\, YX \hbox{ both upper triangular} \big\} $$
which we'll call the {\dfn upper-upper scheme}.

\begin{Proposition}\proplabel{D0insideE}
  Recall the scheme 
  $$ D^0 := 
\lim_{z\to 0} \big\{ (X',Y') : X'Y' = \rhocek(z) Y'X' \rhocek(z^{-1}) 
\hbox{ off the diagonal}
\big\} $$
  defined as the flat limit of this one-parameter family.

  This scheme embeds in $E$, via the map
  $$ \tau: (X,Y) \in D^0 \quad \mapsto \quad (w_0 X, Y w_0) \in E $$
  where $w_0$ is the permutation matrix with $1$s along the antidiagonal.
\end{Proposition}

\begin{proof}
  Consider the equations defining $D^z$. When we conjugate $Y'X'$ by
  $\rhocek(z)$, it multiplies the upper triangle by positive powers of $z$,
  the lower triangle by negative powers (and leaves the diagonal alone,
  though we don't care). So the upper triangle of $X'Y'$ is equal to 
  the upper triangle of $Y'X'$ times positive powers of $t$, and after
  rescaling,  the lower triangle of $Y'X'$ is equal to 
  the lower triangle of $X'Y'$ times positive powers of $z$.
  In the limit, we get the equations $X'Y'$ lower triangular, $Y'X'$ 
  upper triangular on $D^0$. Therefore $(w_0 X', Y' w_0) \in E$.
\end{proof}

We will eventually prove in \thmref{degen} that this map is an isomorphism.
We've twisted this scheme by $w_0$ only because it's less confusing to
deal always with upper triangular matrices, rather than to mix upper and lower.

\subsection{Dimensions of the components.}\sseclabel{degencomponents}

The upper-upper scheme $E$ carries an obvious action of 
pairs of upper triangular matrices, $B_+\times B_+$:
$$ (U_1, U_2) \cdot (X,Y) := (U_1\, X\, U_2^{-1},\, U_2\, Y\, U_1^{-1}) $$
The projection $p: (X,Y) \mapsto X$ is $B_+\!\times\! B_+$-equivariant 
with respect to the action $(U_1, U_2) \cdot X := U_1\, X\, U_2^{-1}$
on the space of single matrices.

We know the orbits of $B_+\times B_+$ on the space of matrices: each orbit
contains a unique {\dfn partial permutation matrix}, a $0,1$-matrix with
at most one $1$ in any row and column. 

We will need also a slightly more specific fact, that each orbit of
$N_+\times N_+$ contains a unique {\dfn monomial matrix}, which has at
most one nonzero entry in each row and column.

Our reference for facts about these
orbits, in particular their dimensions, is \cite{MS}:

\begin{Proposition}[\cite{MS}, theorem 15.28]\proplabel{ppdim}
  Let $\pi$ be a partial permutation matrix. The dimension of
  $B_+ \pi B_+$ inside $M_n(\complexes)$ is the number of matrix entries
  such that a $1$ entry in $\pi$ is either on it, directly below,
  or directly to the left.
\end{Proposition}

Given a partial permutation matrix $\pi$, let
$$ E_\pi := p^{-1}\big( (B_+\times B_+)\cdot \pi\big) $$
so the $\{E_\pi\}$ give a finite decomposition of $E = \coprod_\pi E_\pi$
into $B_+\times B_+$-invariant locally closed subsets.

\begin{Lemma}\lemlabel{Ecomponents}
  The stratum $E_\pi$ is smooth and irreducible, 
  of dimension $n^2 + ($rank of $\pi)$. 

  If $\pi$ is a permutation matrix (not just partial), then the set
  $$ (N_+\times N_+) \cdot \big\{ (\pi t, s \pi^{-1}), \quad 
  s,t\hbox{ invertible diagonal} \big\} $$
  is an open dense subset of $E_\pi$.
\end{Lemma}

\begin{proof}
  The fiber over the ``central'' point $\pi \in E_\pi$ is
  $$ \{ Y: \pi Y, \, Y\pi\hbox{ both upper triangular} \} $$
  which is a vector space. Since $B_+\times B_+$ acts transitively
  on $p(E_\pi)$, $E_\pi$ is a vector bundle over the 
  smooth irreducible $p(E_\pi)$,
  and therefore smooth and irreducible.
  Moreover, the dimension of $E_\pi$ is the dimension of the orbit
  in the base (given by \propref{ppdim}), plus the dimension of the fiber.

  So let's compute the fiber dimension. The conditions $\pi Y, Y\pi$
  upper triangular become, on matrix entries, 
  that $Y_{ij}$ must be zero if there is a $1$ in $\pi$ directly
  (and strictly) to the left, or directly (and strictly) below, entry $ij$. 
  Otherwise $Y_{ij}$ is free, 
  which includes the case when there is a $1$ in $\pi$ actually in entry $ij$.
  
  Every matrix entry therefore ``counts'' for the dimension of the
  base (by \propref{ppdim}) if it has a $1$ in $\pi$ below or to the
  left, counts for the fiber if it doesn't, and counts for both if it
  is actually placed at a $1$ in $\pi$. So the total count is $n^2$ plus 
  the number of $1$s, as was to be shown.

  Consider now the $N_+\times N_+$ orbit of the point $(\pi t, s \pi^{-1})$,
  and assume $s$ invertible, and that $s^{-1} t$ has no repeated entries. 
  (Even excluding those $(s,t)$ it turns out that 
  we'll still get a dense open set.)
  Plainly, this orbit is contained in $E_\pi$ by its definition. 
  The infinitesimal stabilizer of $(\pi t, s \pi^{-1})$ consists of 
  those pairs $(A,B)$ of strictly upper-triangular matrices such that
  $$ A\pi t - \pi t B = 0, \quad    B s \pi^{-1} - s \pi^{-1} A = 0. $$
  So 
  $$ \pi^{-1} A \pi = t B t^{-1} = s^{-1} B s $$
  making $B$ commute with the generic diagonal matrix $s^{-1} t$.
  Therefore $B$ is diagonal, hence zero, and so too is $A$.

  Since the $N_+\times N_+$-stabilizer of $(\pi t, s \pi^{-1})$
  is trivial, its orbit is $2{n\choose 2}$ dimensional. No two of these
  orbits intersect (see the comment before the lemma about monomial
  matrices), so we have a $2n$-dimensional family of them, in all
  comprising $n^2+n$ dimensions. This is the same dimension as $E_\pi$,
  so this subset is open (hence dense, since $E_\pi$ is irreducible).
\end{proof}

The following lemma tells us some (but not all) of the equations
separating the components of $E$.

\begin{Lemma}\lemlabel{recognizeEpi}
  Let $\pi$ be a permutation matrix. Then 
$$(X,Y) \in E_\pi \quad\Longrightarrow\quad \diag(XY) = \pi\cdot \diag(YX), $$
  i.e. $(XY)_{ii} = (YX)_{\pi(i),\pi(i)}$, $i=1\ldots n$.
\end{Lemma}

\begin{proof}
  It's enough to test this equality on the dense subset given us by
  \lemref{Ecomponents}, consisting of elements of the form
  $(X,Y) = (U_1 \pi t U_2^{-1}, U_2 s \pi^{-1} U_1^{-1})$, 
  where $U_1,U_2 \in N_+$ and $s,t$ are diagonal.
  $$ XY = U_1 \pi t s \pi^{-1} U_1^{-1} $$
  $$ YX = U_2 s t U_2^{-1} $$
  So their diagonals are the same as those of $\pi ts \pi^{-1}$ and $st$.
\end{proof}

We give a precise conjecture of the equations defining the closure of
$E_\pi$ in \secref{conjectures}. Note that the obvious component of $E$,
in which both $X$ and $Y$ are themselves 
upper triangular, is $\overline{E_1}$, 
whereas the component that interests us most is $\overline{E_{w_0}}$.

\begin{proof}[Proofs of theorems \ref{thm:degen} and \ref{thm:comps}.]
  The scheme $E$ is defined by $n^2-n$ equations, and by
  \lemref{Ecomponents} is a finite union of pieces $\{E_\pi\}$ each of
  codimension $\geq n^2-n$.  So it is a complete intersection.

  Therefore it is pure, and only those
  pieces $\overline{E_\pi}$ of codimension exactly $n^2-n$ are components 
  (the others lie in the closure).
  These are the $\overline{E_\pi}$ for which $\pi$ has rank $n$, i.e. is a
  permutation matrix and not just a partial permutation matrix.
  In particular this proves the first statement in \thmref{comps}.

  It remains to show that $E$ is reduced. Since it is a complete
  intersection and therefore Cohen-Macaulay, being generically
  reduced implies that it is reduced (see \cite{E}, exercise 18.9).
  We will now find a smooth, reduced point $(\pi t, s \pi^{-1})$ on
  each $E_\pi$.

  Let $t,s$ be generic diagonal matrices (the genericity condition 
  will be specified in due course). 
  The scheme $E$ is the zero set of the map
  $$ (X,Y) \mapsto \hbox{ strict lower triangles of } XY,YX,$$
  whose differential at the point $(\pi t, s \pi^{-1})$ is
  \begin{eqnarray*}
    (A,B) &\mapsto& \hbox{ strict lower triangles of } 
        A Y + X B, Y A + B X \\
        &=& \hbox{ strict lower triangles of } 
        A s\pi^{-1} + \pi t B, \,\, s \pi^{-1} A + B \pi t.
  \end{eqnarray*}
  We want to show this differential is onto. 

  Consider $A = \lambda e_{\pi(i),j}$, $B = \mu e_{k,\pi(l)}$. Then
  \begin{eqnarray*}
    (A,B) &\mapsto& \hbox{ strict lower triangles of } 
        \lambda e_{\pi(i),j} s \pi^{-1} + \mu \pi t e_{k,\pi(l)}, \,
        s \pi^{-1} \lambda e_{\pi(i),j} + \mu e_{k,\pi(l)} \pi t \\
    &=& \hbox{ strict lower triangles of } 
        \lambda s_j e_{\pi(i),\pi(j)} + \mu t_k e_{\pi(k),\pi(l)}, \, 
        \lambda s_i e_{ij} + \mu t_l e_{kl} 
  \end{eqnarray*}
  In particular $\lambda=1, \mu=0$ gives us
  $$
  (A,B) \mapsto \hbox{ strict lower triangles of } 
         s_j e_{\pi(i),\pi(j)}, s_i e_{ij}.
  $$
  If $i>j$ but $\pi(i)<\pi(j)$, we can use this to produce pairs
  $(0, e_{ij})$. 
  If $i<j$ but $\pi(i)>\pi(j)$, we can use this to produce pairs
  $(e_{\pi(i),\pi(j)},0)$.

  The hard case is when  $i>j$ and $\pi(i)>\pi(j)$, then $(i,j)=(k,l)$ gives us
  $$
  (A,B) \mapsto \hbox{ strict lower triangles of } 
        (\lambda s_j + \mu t_i) e_{\pi(i),\pi(j)},
        (\lambda s_i + \mu t_j) e_{ij} 
  $$
  and as long as $s_j/s_i \neq t_i/t_j$ for any $i,j$, we can adjust
  $\lambda,\mu$ to obtain $(0, e_{ij})$ and $(e_{\pi(i),\pi(j)}, 0)$
  in the image. So we've gotten every pair of matrices where one has
  zero strict lower triangle and the other has exactly one entry in 
  the strict lower triangle. These generate the target so the differential
  is indeed onto.

  We've found a reduced point in each component $E_\pi$ of $E$. 
  Therefore $E$ is generically reduced, hence by its Cohen-Macaulayness
  it's reduced.
  
  We're now ready to knock off \thmref{degen} and the remainder of
  \thmref{comps}.  \Thmref{degen} amounts to the statement that the
  map in \propref{D0insideE} is an isomorphism, which we'll now prove.

  Since $D^1$ and $E$ are complete intersections defined by $n^2-n$
  quadratics, they have the same degree, $2^{n^2-n}$.
  This map $\tau: (X,Y) \mapsto (w_0 X, Y w_0)$ from
  \propref{D0insideE} is linear, so preserves degree. Taking flat
  limits also preserves degree.  So the image $\tau(D^0)$ inside $E$ has
  the same degree as $E$.  Since $E$ is equidimensional of the same
  dimension as $\tau(D^0)$, we find that $\tau(D^0)$ must include all of
  $E$'s components. But this lower bound on $\tau(D^0)$ is already $E$,
  since $E$ is reduced.

  By \cite{KS}, the (reduction of the) $z\to 0$ limit of any
  component of $D^z$ is again equidimensional, hence a union of some
  components of $E$. We want to see that the commuting component of $D$
  limits only to $\overline{E_1}$, and that the non-commuting components
  of $D$ accounts for all the other components of $E$.
  A priori one might expect some components $\{\overline E_\pi\}$ to arise
  as components of both limits, but as $E$ is generically reduced
  this does not happen.

  Let $t,s$ be generic and $\pi\neq 1$. 
  Now note that each point $(\pi t, s \pi^{-1})$ is in every $D^z$. 
  For $z=1$, this point is in a 
  non-commuting component of $D = D^1$. For $z=0$, this point is in the
  $\overline{E_\pi}$ component (and no other, by the genericity). 
  So the limit of the non-commuting components of $D$ is all the
  non-identity components of $E$.
\end{proof}

Unfortunately, the result of \cite{KS} doesn't let us rule out the
possibility that the $z\to 0$ limit of the commuting variety has
embedded components. To show this doesn't happen, it would be enough
to know that the variety $\overline{E_{w_0}}$ is defined by the additional 
equations $\diag(XY) = w_0\cdot \diag(YX)$ with no more needed. 
That would also imply that the commuting scheme is reduced, which
remains unknown at the time of this writing.

\subsection{Some results about (multi)degrees of components.}%
\sseclabel{numerology}

Let $d_\pi$ denote the degree of the homogeneous
affine variety $\overline{E_\pi}$.

\begin{Proposition}\proplabel{degrees}
  \begin{itemize}
  \item The sum $\sum_{\pi\in S_n} d_\pi$ is $2^{n^2-n}$.
  \item Denote by $\star: S_k \times S_{n-k} \to S_n$ the standard
    concatenation of permutations.\\
    Then $d_{\pi\,\star\,\rho} = d_\pi d_\rho$.
  \item If $w_0$ is the permutation of length $n$ of maximum length, then
    $$ d_\pi = d_\pi^{-1} = d_{w_0 \pi w_0} = d_{w_0 \pi^{-1} w_0}.$$
  \end{itemize}
\end{Proposition}

\begin{proof}
  For the first statement, the right-hand side is the degree of
  the quadratic complete intersection $E$, which is the sum of the
  degrees of its components.

  For the second, note that 
  \begin{eqnarray*}
\overline{E_\pi} \times\overline{E_\rho} 
        \times ({\mathbb A}^{k\times (n-k)})^2
  &\to& \overline{E_{\pi\star\rho}} \\
\big( (X_1,Y_1), (X_2,Y_2), M_1, M_2) 
  &\mapsto& 
  \left( 
    \left(\begin{array}{cc}X_1 & M_1 \\
                            0  & Y_1\end{array}\right),
    \left(\begin{array}{cc}X_2 & M_2 \\
                            0  & Y_1\end{array}\right) 
  \right) 
\end{eqnarray*}
   is an isomorphism, and linear so degree-preserving. 

   The third is really two statements. The $\pi\leftrightarrow \pi^{-1}$
   symmetry comes from the map $(X,Y) \mapsto (Y,X)$.
   The map
   $$ (X,Y) \mapsto (w_0 X^T w_0, w_0 Y^T w_0) $$
   is also easily seen to give a linear isomorphism of $E_\pi$ and
   $E_{w_0 \pi^{-1} w_0}$.
\end{proof}

With Macaulay 2 \cite{M2}, we computed these degrees for small $n$ and
$\star$-irreducible $\pi$:
$$ d_{1} = 1 $$
$$ d_{21} = 3 $$
$$ d_{231} = d_{312} = 13, \quad d_{321} = 31 $$
So for example when $n=3$, 
\begin{eqnarray*}
 2^{3^2 - 3} &=& d_{123} + d_{213} + d_{132} + d_{312} + d_{231} + d_{321} \\
&=& d_1^3 + d_{21} d_1 + d_1 d_{21} + d_{312} + d_{231} + d_{321} \\
&=& 1 + 3 + 3 + 13 + 13 + 31.
\end{eqnarray*}


In fact one can sharpen \propref{degrees} a great deal using the theory of
multidegrees \cite{MS} (also known as equivariant multiplicities \cite{Ro}),
using not only the rescaling action but the full torus action on $D^0$.
Since we don't want to recapitulate this theory here -- except to say that
it assigns each cycle a homogeneous polynomial, rather than just a number --
we give only a little taste, using the $2$-torus action that scales
$X$ and $Y$ individually. 

\begin{Proposition}\proplabel{multidegrees}
  Let $A,B$ be the usual generators of
  the weight lattice of the $2$-torus scaling $X$,$Y$ individually.
  Let $d'_\pi$ denote the bidegree of $E_\pi$, a homogeneous 
  polynomial in $\naturals[A,B]$.
  \begin{itemize}
  \item The sum $\sum_{\pi\in S_n} d'_\pi$ is $(A+B)^{n^2-n}$.
  \item Denote by $\star: S_k \times S_{n-k} \to S_n$ the standard
    concatenation of permutations.  Then 
    $$ d'_{\pi\star\rho} = d'_\pi\, d'_\rho\, (AB)^{k(n-k)}. $$
  \item If $w_0$ is the permutation of length $n$ of maximum length,
    then $d'_\pi = d'_{w_0 \pi^{-1} w_0}$.
  \item $d'_\pi(A,B) = d'_{\pi^{-1}}(B,A).$
  \end{itemize}
\end{Proposition}
The proofs are exactly as in \propref{degrees}.
The $n=3$ example now becomes
$$ d'_{1} = 1,\quad d'_{21} = A^2 + AB + B^2 $$
$$
d'_{231} = 2 A^5 B + 4 A^4 B^2 + 4 A^3 B^3 + 2 A^2 B^4 + A B^5, \quad
d'_{312} = A^5 B + 2 A^4 B^2 + 4 A^3 B^3 + 4 A^2 B^4 + 2 A B^5 $$
$$ d'_{321} =A^6 + 3 A^5 B + 7 A^4 B^2 + 9 A^3 B^3 + 7 A^2 B^4 + 3 A B^5 + B^6
$$

\begin{eqnarray*}
 (A+B)^{3^2 - 3} 
&=& d'_{123} + d'_{213} + d'_{132} + d'_{312} + d'_{231} + d'_{321} \\
&=&(AB)^3 (d'_1)^3 + AB d'_{21} d'_1 + AB d'_1 d'_{21} 
        + d'_{312} + d'_{231} + d'_{321}\\
&=& A^3 B^3 + 2 (A^3 B + A^2 B^2 + A B^3) \\
&+& (2 A^5 B + 4 A^4 B^2 + 4 A^3 B^3 + 2 A^2 B^4 + A B^5) \\
&+& (A^5 B + 2 A^4 B^2 + 4 A^3 B^3 + 4 A^2 B^4 + 2 A B^5) \\
&+& (A^6 + 3 A^5 B + 7 A^4 B^2 + 9 A^3 B^3 + 7 A^2 B^4 + 3 A B^5 + B^6)
\end{eqnarray*}

\section{Conjectures}\seclabel{conjectures}

There are two main conjectures about the commuting scheme: 
that it is reduced, and that it is Cohen-Macaulay. We state some
conjectures sharpening these two.

\begin{Conjecture}
  The variety $\overline{E_\pi}$ of $E$ is defined as a scheme 
  by three sets of equations:
  \begin{enumerate}
  \item those defining $E$, which say $XY,YX$ upper triangular
  \item those given by \lemref{recognizeEpi}, 
    that $\diag(XY) = \pi\cdot \diag(YX)$
  \item those defining the $\pi,\pi^{-1}$ matrix Schubert varieties:
    for each pair $i,j$
    the rank of the lower-left $i\times j$ rectangle in $X$ (resp. in $Y$)
    is bounded above by the number of $1$s in that rectangle in $\pi$
    (resp. in $\pi^{-1}$).
  \end{enumerate}
  Note that for $\pi=w_0$, the third set is empty.

  Moreover, if we impose just the first and third set of equations,
  we get the reduced scheme $\overline{ \union_{\rho\leq \pi} E_\rho}$.
\end{Conjecture}

If this is proved for $\pi=w_0$, it implies that the commuting scheme
(to which $\overline E_{w_0}$ deforms) is reduced, i.e. is the
commuting variety.

\begin{Conjecture}
  Each individual $\{E_\pi\}$, and each union
  $\overline{ \union_{\rho\leq \pi} E_\rho}$ , is Cohen-Macaulay.
\end{Conjecture}

Note that these statements are trivial for the component $\overline{E_1}$, 
being a linear subspace $\{(X,Y) :$ both upper triangular$\}$.
Perhaps they can be proved by induction in the Bruhat order.

We repeat our earlier-stated conjecture (which doesn't seem
to imply anything directly about the commuting scheme):

\begin{Conjecture}
  For $\lie{g}$ a reductive Lie algebra, the diagonal commutator scheme
  of $\lie{g}$ is a reduced complete intersection with two components.
\end{Conjecture}

In the $\gln$ case, Terry Tao conjectured in a discussion the equations
defining the other component. First, we find some equations that 
do in fact hold.

\begin{Proposition}\proplabel{othercomp}
  Consider the $n\times 2n$ matrix, whose first $n$ columns are the
  diagonals of $X^i$, $i=0\ldots n-1$, and next $n$ are the diagonals
  of $Y^i$, $i=0\ldots n-1$. 

  If $(X,Y) \in H$ but $[X,Y] \neq 0$, then the rank of this $n\times 2n$
  matrix is at most $n-1$. In particular every size $n$ minor vanishes. 
\end{Proposition}

\begin{proof}
  Let $K = [X,Y]$. Then the nonzero diagonal matrix $K$ 
  is trace-perpendicular to any element $Z_X$ in the centralizer $C_X$ of $X$:
  $$ \Tr (K Z_X) = \Tr ([X,Y] Z_X) = \Tr ([Z_X,X] Y) = \Tr 0 = 0 $$
  The same argument holds for any $Z_Y$ in the centralizer $C_Y$ of
  $Y$ (rotating the opposite direction), and any linear combination
  $Z_X + Z_Y$.

  The functional $\Tr(K\cdot)$ is only sensitive to the diagonal entries,
  and the trace form $\Tr(\cdot\cdot)$ is nondegenerate. So the
  projection ``take diagonals'' from $C_X + C_Y$ to $\lie{t}$ is not onto,
  since it only hits $K^\perp$. (This is where we use $K\neq 0$.)

  The argument so far would work fine in any semisimple $\lie{g}$, with the
  Killing form in place of the trace form. In the $\gln$ case, we have a
  bunch of matrices we know to be in $C_X$ (resp. $C_Y$), 
  namely the powers of $X$ (resp. $Y$). The non-ontoness of the projection
  then gives us the rank claim in the proposition.
\end{proof}

Note that this gives one equation each on $X$ and $Y$ individually -- while
every $X$ commutes with some $Y$ (e.g. $0$ or $X$ itself), not every $X$
has a {\em nonzero} diagonal commutator with some $Y$.

\begin{Conjecture}
  The equations in \propref{othercomp} define the other component(s)
  of the diagonal commutator scheme.
\end{Conjecture}

With Macaulay 2, we verified this in the $\lie{gl}_3$ case -- first
by finding the equations, then using them to suggest the conjecture.
The other conjectures were also all verified in this case.

\bibliographystyle{alpha}    

\end{document}